\documentclass{amsart}
\usepackage{amsfonts,amscd}

\newtheorem{claim}{}[section]
\newtheorem{theorem}[claim]{Theorem}
\newtheorem{lemma}[claim]{Lemma}
\newtheorem{proposition}[claim]{Proposition}
\newtheorem{corollary}[claim]{Corollary}
\newtheorem{remark}[claim]{Remark}

\newtheorem{example}[claim]{Example}

\def\proclaim #1. #2\par{\medbreak
\noindent{\bf#1.\enspace}{\sl#2}\par\medbreak}
\makeatother

\DeclareMathOperator{\spn}{span}
\DeclareMathOperator{\wk}{wk}

\DeclareMathOperator{\WOT}{WOT}

\DeclareMathOperator{\Real}{Re}
\DeclareMathOperator{\ran}{ran} 
\DeclareMathOperator{\ext}{ext}
\DeclareMathOperator{\ball}{ball}

\begin{document}

\title{One-Sided Projections on $C^*$-Algebras}

\date{Mar. 5, 2002}


\author{David P. Blecher}
\address{Department of Mathematics, University of Houston, Houston, TX 
77204-3008}
\email[David P. Blecher]{dblecher@math.uh.edu}

%\author{Edward G. Effros}
%\address{Department of Mathematics, UCLA, Los Angeles, CA 90095-1555}
%\email[Edward G. Effros]{ege@math.ucla.edu}

\author{Roger R. Smith}
\address{Department of Mathematics, Texas A\&M University, College Station, TX 77843-3368}
\email[Roger R. Smith]{rsmith@math.tamu.edu}

\author{Vrej Zarikian}
\address{Department of Mathematics, University of Texas, Austin, TX 78712-1082}
\email[Vrej Zarikian]{zarikian@math.utexas.edu}

\thanks{Blecher and Smith were partially supported by the National Science Foundation. Zarikian was partially supported by a UCLA Dissertation Year Fellowship.}

\begin{abstract}
In \cite{BEZ}, the notion of a complete one-sided $M$-ideal for an operator space $X$ was introduced as a generalization of Alfsen and Effros' notion of an $M$-ideal for a Banach space \cite{AE72}. In particular, various equivalent formulations of complete one-sided 
$M$-projections were given. In this paper, some sharper equivalent formulations are given in the special situation that $X = \mathcal{A}$, a $C^*$-algebra 
(in which case the complete left $M$-projections are simply left multiplication on $\mathcal{A}$ by a fixed orthogonal projection in $\mathcal{A}$ or its multiplier algebra). The proof of the first equivalence makes use of a technique which is of interest in its own right---a way of ``solving'' multi-linear equations in von Neumann algebras. This technique is also applied to show that preduals of von Neumann algebras have no nontrivial complete one-sided $M$-ideals. In addition, we show that in a $C^*$-algebra, the intersection of finitely many complete one-sided $M$-summands need not be a complete one-sided $M$-summand, unlike the classical situation.
\end{abstract}
\maketitle

%%%%%%%%%%%%%%%%%%%%%%%%%%%%%%%%%%%%%%%%%%%%%%%%%%%%%%%%%%%%%%%%%%%%%%

\let\text=\mbox

\section{Introduction}

One of the main results of the present paper is the following:
 
\begin{theorem}  \label{main}
Let $\mathcal{A}$ be a $C^*$-algebra, and let $P : \mathcal{A} \rightarrow
\mathcal{A}$ be a idempotent linear map.  The following are equivalent:
\begin{itemize}
\item [(i)]  there is an orthogonal projection $e \in M(\mathcal{A})$ 
such that $Px = ex$  for all $x \in \mathcal{A}$,
\item [(ii)] the map $x \mapsto \begin{bmatrix} Px \\ x - Px \end{bmatrix}$ is isometric as a map from $\mathcal{A}$ to $C_2(\mathcal{A})$,
\item [(iii)] the map $\begin{bmatrix} x \\ y \end{bmatrix}
\mapsto \begin{bmatrix} Px \\ y \end{bmatrix}$ on $C_2(\mathcal{A})$  is contractive,
\item [(iv)] the map $(Px,y-Py) \mapsto \begin{bmatrix} Px & y-Py \end{bmatrix}$ is isometric as a map from \linebreak $P\mathcal{A} \oplus_\infty (Id - P)\mathcal{A}$ to $R_2(\mathcal{A})$,
\item [(v)] $(Px)^* Px \leq x^* x$ for all $x \in \mathcal{A}$,
\item [(vi)] $(Px)^*(y - Py) = 0$ for all $x, y \in \mathcal{A}$.
\end{itemize}
\end{theorem}

In the above, and throughout this paper, $C_2(\mathcal{A})$ and $R_2(\mathcal{A})$ are the first column and row of $M_2(\mathcal{A})$, respectively. Since $M_2(\mathcal{A})$ is a $C^*$-algebra, $C_2(\mathcal{A})$ and $R_2(\mathcal{A})$ have canonical norms. Indeed, the norm of a column (resp. row) with entries $x$ and $y$ from $\mathcal{A}$ is $\sqrt{\|x^*x + y^*y\|}$ (resp. $\sqrt{\|xx^* + yy^*\|}$). We write $M(\mathcal{A})$ for the multiplier algebra of $\mathcal{A}$ (see \cite{Ped} for example). Of course $M(\mathcal{A}) = \mathcal{A}$ if $\mathcal{A}$ has an identity.
An {\em idempotent} map is a map $P$ for which $P \circ P = P$. A {\em projection} on a Hilbert space, or in a $C^*$-algebra, will mean an orthogonal projection. If $e$ is such a projection, then $e^\perp$ will denote the complementary projection, $1 - e$. Throughout we use the symbol $\mathfrak{H}$ for a Hilbert space. 

Let $X$ be a general operator space (i.e. a closed linear subspace of 
some $B(\mathfrak{H})$). A linear map $T$ from $X$ to another operator space $Y$ induces a canonical map $T_n:M_n(X) \to M_n(Y)$. We say that $T$ is {\em $n$-isometric} if $T_n$ is an isometry, and that $T$ is {\em completely isometric} if it is $n$-isometric for every $n \in \mathbb{N}$. In \cite{BEZ}, we defined a {\em complete left $M$-projection} on $X$ (resp. {\em left $M$-projection} on $X$) to be an idempotent linear map $P$ on $X$ such that the associated map in (ii) above is {\em completely isometric} (resp. isometric). Indeed, in the rest of the paper we will write $\nu_P^c$ for this map (from $X$ to $C_2(X)$). We showed in  \cite{BEZ} that the complete left $M$-projections on a $C^*$-algebra $\mathcal{A}$ are exactly the maps in (i) of the Theorem above. Thus the equivalence of (i) and (ii) may be interpreted as the statement that every left $M$-projection on a $C^*$-algebra is a complete left $M$-projection. Note also that (ii) is a non-commutative analog of the formula
\[
       \|x\| = \max\{\|Px\|,\|x - Px\|\},
\]
which characterizes classical $M$-projections on a Banach space \cite{AE72,HWW93}.

In \cite{BEZ} we proved that $P$ being a
complete left $M$-projection on a general 
operator space $X$ is equivalent to a number of 
other conditions.  For example, it is equivalent to saying that
the associated map in (iii) of the Theorem above is  
completely contractive.   It is also equivalent to a 
`matricial version' of (v) holding on $X$ (see
\cite{Ble_Shilov}, Section 4).    Condition (iv) in fact 
also characterizes complete left $M$-projections on 
general (possibly non-self-adjoint) operator algebras which
possess a contractive approximate identity, as we shall
see in Section 6.

The other main result in our paper is the fact that the predual or dual of
a von Neumann algebra possesses no nontrivial complete one-sided $M$-projections (or complete one-sided $M$-ideals).

In the final section of the paper we show that the
intersection of finitely many complete one-sided $M$-summands in 
a $C^*$-algebra need not be
a complete one-sided $M$-summand, unlike the classical situation.

For the history of some of the topics discussed in our
paper, and references to some related work of others,
we refer the reader to the discussion in the
introduction of \cite{BEZ}.   One part of this
history should be amplified.    Part (iii) of
our Theorem \ref{main} above, and its generalization below in Lemma \ref{vonNeum},
may be viewed as a strengthening
of the  characterization of `left multipliers of
operator spaces' given in Theorem 4.6 of
\cite{BEZ}, in the particular case that
the operator space is a $C^*$-algebra.  We made much
use of this Theorem 4.6 in \cite{BEZ} and elsewhere.  
Historically, left multipliers for an operator system $S$
were first considered
by W. Werner in \cite{WWerner} around 1998, where he
obtained the version appropriate to operator systems of our
Theorem 4.6 of
\cite{BEZ};  he also characterized left multipliers in terms
of the injective envelope $I(S)$. A year later
the first author, unaware of this work,  considered left multipliers
of operator spaces in \cite{Ble_Shilov} and \cite{BP01}.
Theorem 4.6 in \cite{BEZ} was inspired by W. Werner's
characterization, but it was only recently that we
saw that our result can be deduced from W. Werner's original
theorem from \cite{WWerner}.   Briefly, the point is that any
operator space $X$ may be studied in terms of the off-diagonal
$\begin{bmatrix} 0 & X\\ X^* & 0 \end{bmatrix}$ inside  the `Paulsen system'; and thus may (after making a few appropriate choices) be seen
to fall within the framework of the `non-unital operator systems' that Werner considers in his theorem.

\section{The Glimm-Halpern Reduction Theory} \label{Section:Glimm-Halpern}

In this section we give a brief account of the Glimm-Halpern 
reduction theory (cf. \cite{Gli60} and \cite{Hal69}). 
The use of this theory in the proof of the main result in the next section was generously suggested by Edward Effros, who noticed its effective use in \cite{CS93}.\\

Let $\mathcal{R}$ be a von Neumann algebra. Let $\Omega$ be the spectrum of the center $\mathcal{Z(R)}$. Then 
$\Omega$ is an extremely disconnected compact Hausdorff space (\cite{KR_I}, Theorem 5.2.1), meaning that the closure of every open set is open. We have that $\mathcal{Z(R)} \cong C(\Omega)$ via the Gelfand transform. For each $\omega \in \Omega$, let $\mathcal{M}_\omega = \ker(\omega) \subset \mathcal{Z(R)}$ be the corresponding maximal ideal. Let $\mathcal{I}_\omega \subset \mathcal{R}$ be the norm-closed, 2-sided ideal generated by $\mathcal{M}_\omega$. It is easy to check that
\[
	\mathcal{I}_\omega = \overline{\spn\{zr: z \in \mathcal{M}_\omega, 
		r \in \mathcal{R}\}}.
\]
Define $\mathcal{R}_\omega = \mathcal{R}/\mathcal{I}_\omega$, a $C^*$-algebra. We regard $\{\mathcal{R}_\omega: \omega \in \Omega\}$ as a decomposition of $\mathcal{R}$. It has several attractive features, which we now discuss. We include the simpler proofs.\\

First some notation: for each $x \in \mathcal{R}$ and each $\omega \in \Omega$, let $x(\omega) = x + \mathcal{I}_\omega \in \mathcal{R}_\omega$.

\begin{theorem}[\cite{Gli60}, Remarks before Lemma 9] \label{Glimm_thm_1}
Let $x \in \mathcal{R}$. Then $\|x\| = \sup\{\|x(\omega)\|: \omega \in \Omega\}$.
\end{theorem}

\begin{proof}
Clearly, $\sup\{\|x(\omega)\|: \omega \in \Omega\} \leq \|x\|$. On the other hand, there exists a pure state $\varphi:\mathcal{R} \to \mathbb{C}$ such that $\|x\| = \varphi(x^*x)^{1/2}$ (\cite{KR_I}, Theorem 4.3.8). We claim that $\omega = \varphi|_{\mathcal{Z(R)}} \in \Omega$. Let us assume the claim for the moment. Because $\varphi|_{\mathcal{M}_\omega} = \omega|_{\mathcal{M}_\omega} = 0$, 
it follows that 
$\varphi|_{\mathcal{I}_\omega} = 0$, by the Cauchy-Schwarz inequality. 

Thus for any $y \in \mathcal{I}_\omega$,
\[
	\|x + y\|^2 \geq \varphi((x + y)^*(x + y)) = \varphi(x^*x) = \|x\|^2,
 \]
giving
 $\|x(\omega)\| \geq \|x\|$.

We now prove the claim. Let $(\pi_\varphi,\mathfrak{H}_\varphi)$ be the GNS construction for $\mathcal{R}$ corresponding to $\varphi$. Because $\varphi$ is pure, $\pi_\varphi(\mathcal{R}) \subset B(\mathfrak{H}_\varphi)$ is irreducible (\cite{KR_II}, Theorem 10.2.3). 
Thus $\mathbb{C} \subset \pi_\varphi(\mathcal{Z(R)}) \subset \mathcal{Z}(\pi_\varphi(\mathcal{R})) = \mathbb{C}$ (\cite{KR_I}, Theorem 5.4.1),
which implies that 
$\mathcal{Z(R)}/(\ker(\pi_\varphi) \cap \mathcal{Z(R)}) \cong \mathbb{C}$. Now let $z_1, z_2 \in \mathcal{Z(R)}$. By the previous discussion, there exist $\lambda_1, \lambda_2 \in \mathbb{C}$ and $u_1, u_2 \in \ker(\pi_\varphi) \cap \mathcal{Z(R)}$ such that $z_1 = \lambda_1 + u_1$ and $z_2 = \lambda_2 + u_2$. Since $\ker(\pi_\varphi) \subset \ker(\varphi)$,
\[
	\varphi(z_1z_2) = \lambda_1\lambda_2 = \varphi(z_1)\varphi(z_2).
\]
\end{proof}

\begin{corollary} \label{Glimm_cor_1}
Let $x, y \in \mathcal{R}$. Then $x = y$ if and only if $x(\omega) = y(\omega)$ for all $\omega \in \Omega$.
\end{corollary}

\begin{corollary} \label{Glimm_cor_2}
Let $x \in \mathcal{R}$. Then $x \in \mathcal{Z(R)}$ if and only if $x(\omega) \in \mathbb{C}$ for all $\omega \in \Omega$.
\end{corollary}

\begin{proof}
Suppose $x \in \mathcal{Z(R)}$. Let $\omega \in \Omega$ be arbitrary. Then $\omega(x - \omega(x)) = 0$ $\Rightarrow$ $x - \omega(x) \in \mathcal{M}_\omega \subset \mathcal{I}_\omega$. Thus, $x(\omega) = \omega(x)$. Conversely, suppose $x(\omega)$ is a scalar for all $\omega \in \Omega$. Then for any $y \in \mathcal{R}$,
\[
	(xy)(\omega) = x(\omega)y(\omega) = y(\omega)x(\omega) = (yx)(\omega)
\]
for all $\omega \in \Omega$.  Thus $xy = yx$.
\end{proof}

\begin{corollary} \label{Glimm_cor_3}
Let $f \in C(\Omega)$. Then there exists an $x \in \mathcal{Z(R)}$ such that $x(\omega) = f(\omega)$ for all $\omega \in \Omega$.
\end{corollary}

\begin{proof}
Recall that $\mathcal{Z(R)} \cong C(\Omega)$ via the Gelfand transform. Thus, there exists an $x \in \mathcal{Z(R)}$ such that $\omega(x) = f(\omega)$ for all $\omega \in \Omega$. By the proof of the previous corollary, $x(\omega) = \omega(x)$ for all $\omega \in \Omega$.
\end{proof}

\begin{theorem}[\cite{Gli60}, Lemma 10] \label{Glimm_thm_2}
Let $x \in \mathcal{R}$. Then the map $\Omega \to \mathbb{C}:\omega \mapsto \|x(\omega)\|$ is continuous.
\end{theorem}

\begin{theorem}[\cite{Hal69}, Theorem 4.7] \label{Halpern_thm_1}
Let $\omega \in \Omega$. Then $\mathcal{I}_\omega$ is primitive (i.e. $\mathcal{R}_\omega$ has a faithful irreducible representation).
\end{theorem}

\section{A Theorem on ``Solving'' Multi-linear 
Equations in von Neumann Algebras} \label{Section:Smith}

In recent years, several authors have considered `elementary operators' on 
von Neumann algebras (see e.g. \cite{Mathieu}), and their relation
to multi-linear equations in von Neumann algebras.  In 
this section we state and prove a theorem on ``solving'' 
such equations. 
Many of the preliminary steps used in 
reaching this main result are themselves quite interesting.
   
Let $\mathcal{A}$ be a $C^*$-algebra. Then the map $\theta:\mathcal{A} \otimes \mathcal{A} \to CB(\mathcal{A})$ defined by
\[
	\theta(\sum_{j=1}^n a_j \otimes b_j)(x) = \sum_{j=1}^n a_jxb_j
\]
for all $x \in \mathcal{A}$ is well-defined.
A result of Chatterjee and Smith (\cite{CS93}, 
Lemma 2.2) identifies the kernel of an extension of
$\theta$ to a
certain completed tensor product,
in the case that $\mathcal{A} = \mathcal{R}$, a von Neumann algebra. 
We shall only need the following algebraic result, for which a 
streamlined proof is provided.

\begin{lemma}[Chatterjee-Smith] \label{Chatterjee-Smith}
Let $\mathcal{R}$ be a von Neumann algebra and $\theta:\mathcal{R} \otimes 
\mathcal{R} \to CB(\mathcal{R})$ be as in the preceding discussion. Let
\[
	\mathcal{J} = \spn\{az \otimes b - a \otimes zb: a, b \in \mathcal{R},
		z \in \mathcal{Z(R)}\}.
\]
Then $\ker(\theta) = \mathcal{J}$.
\end{lemma}

\begin{proof}
Certainly $\mathcal{J} \subset \ker(\theta)$. On the other hand, suppose $u = \sum_{k=1}^n r_k \otimes s_k \in \mathcal{R} \otimes \mathcal{R}$ and $\theta(u) = 0$. Letting
\[
	r = 
	\begin{bmatrix}
		r_1 & r_2 & \dots & r_n\\
		0 & 0 & \dots & 0\\
		\vdots & \vdots & & \vdots\\
		0 & 0 & \dots & 0
	\end{bmatrix}
	\text{ and }
	s =
	\begin{bmatrix}
		s_1 & 0 & \dots & 0\\
		s_2 & 0 & \dots & 0\\
		\vdots & \vdots & & \vdots\\
		s_n & 0 & \dots & 0
	\end{bmatrix},
\]
one has that $r(x \otimes I_n)s = 0$ for all $x \in \mathcal{R}$. Here $I_n$ is the multiplicative identity of $M_n$. Let $\mathcal{L} = \overline{\spn\{M_n(\mathcal{R})r(\mathcal{R} \otimes I_n)\}}^{\WOT}$, a WOT-closed left ideal of $M_n(\mathcal{R})$. Then (cf. \cite{KR_II}, Theorem 6.8.8) there exists a 
projection $z \in M_n(\mathcal{R})$ such that $\mathcal{L} = M_n(\mathcal{R})z$. In particular, $z \in \mathcal{L}$. We claim that $z \in M_n(\mathcal{Z(R)})$. Indeed, for all $x \in \mathcal{R}$, $\mathcal{L}(x \otimes I_n) \subset \mathcal{L}$, which implies that $z(x \otimes I_n)z = z(x \otimes I_n)$. Replacing $x$ by $x^*$, we conclude that $z(x \otimes I_n) = (x \otimes I_n)z$. Since the choice of $x$ was arbitrary, $z \in (\mathcal{R} \otimes I_n)' = M_n(\mathcal{R}')$ and so $z \in M_n(\mathcal{R}) \cap M_n(\mathcal{R}') = M_n(\mathcal{Z(R)})$, proving the claim. Since $r \in \mathcal{L}$, $rz = r$. Since $\mathcal{L}s = 0$, $zs = 0$. Letting $\odot$ stand for the matrix inner product (cf. \cite{ER_book}, \S 9.1), if $a, b \in M_n(\mathcal{R})$, then
\[
	az \odot b \equiv a \odot zb \mod M_n(\mathcal{J})
\]
since the $(i,j)$ entry of the left-hand side is
\[
	\sum_{k=1}^n (az)_{i,k} \otimes b_{k,j} = \sum_{k=1}^n (\sum_{l=1}^n
		a_{i,l}z_{l,k}) \otimes b_{k,j} = \sum_{k=1}^n \sum_{l=1}^n
		a_{i,l}z_{l,k} \otimes b_{k,j}
\]
and the $(i,j)$ entry of the right-hand side is
\[
	\sum_{l=1}^n a_{i,l} \otimes (zb)_{l,j} = \sum_{l=1}^n a_{i,l} \otimes
		(\sum_{k=1}^n z_{l,k}b_{k,j}) = \sum_{k=1}^n \sum_{l=1}^n
		a_{i,l} \otimes z_{l,k} b_{k,j}.
\]
In particular, letting $0_n$ denote the additive identity of $M_n(\mathcal{R} \otimes \mathcal{R})$, we have that
\[
	u \oplus 0_{n-1} = r \odot s = rz \odot z^\perp s \equiv r \odot 
		z z^\perp s = 0_n \mod M_n(\mathcal{J}).
\]
Hence, $u \in \mathcal{J}$.
\end{proof}

\begin{lemma} \label{vector_sp_lem}
Let $V$ be a vector space and $a, b, c, d, e \in V$, with $e \neq 0$. If
\[
	a \otimes e + e \otimes b + c \otimes d = 0,
\]
then either $a$ and $c$ are multiples of $e$, or $b$ and $d$ are multiples of $e$.
\end{lemma}

\begin{proof}
Let $\varphi$ be a linear functional on $V$ such that $\varphi(e) \neq 0$. Applying $\varphi \otimes Id$ to the above equation yields
\[
	\varphi(a)e + \varphi(e)b + \varphi(c)d = 0.
\]
Thus $b$, $d$, and $e$ are linearly dependent. Let $E = \spn\{b, d, e\}$. Then $\dim(E) = 1 \text{ or } 2$. Case 1: Suppose $\dim(E) = 1$. Then $b$ and $d$ are multiples of $e$. Case 2: Suppose $b$ and $e$ are linearly independent. Then $d = \lambda b + \mu e$. Thus,
\[
	a \otimes e + e \otimes b + c \otimes (\lambda b + \mu e) = 0
\]
$\Rightarrow$
\[
	(a + \mu c) \otimes e + (e + \lambda c) \otimes b = 0
\]
$\Rightarrow$ $c = -\frac{1}{\lambda}e$ and $a = \frac{\mu}{\lambda}e$. Case 3: Suppose $d$ and $e$ are linearly independent. Then $b = \lambda d + \mu e$. Thus
\[
	a \otimes e + e \otimes (\lambda d + \mu e) + c \otimes d = 0
\]
$\Rightarrow$
\[
	(a + \mu e) \otimes e + (\lambda e + c) \otimes d = 0
\]
$\Rightarrow$ $a = -\mu e$ and $c = -\lambda e$.
\end{proof}

\begin{theorem} \label{Smith}
Let $\mathcal{R}$ be a von Neumann algebra and $a, b, c, d \in \mathcal{R}$. If
\[
	ax + xb + cxd = 0
\]
for all $x \in \mathcal{R}$, then there exists a central projection $p \in \mathcal{R}$ such that $pa, pc, p^\perp b, p^\perp d \in \mathcal{Z(R)}$. Conversely, if $pa, pc, p^\perp b, p^\perp d \in \mathcal{Z(R)}$ for a central projection $p \in \mathcal{R}$,
and if $a + b + cd = 0$, then
$ax + xb + cxd = 0$ for all $x \in \mathcal{R}$.
\end{theorem}

\begin{proof}
For the converse, we observe that if 
$pa, pc, p^\perp b, p^\perp d \in \mathcal{Z(R)}$,
then 
\begin{eqnarray*}
ax + xb + cxd &=& p^\perp a x  + x a p 
+ x b p + p^\perp b x + x c d p + p^\perp c x d \\
&=&  p^\perp a x  + x a p
+ x b p + p^\perp b x + x c d p + p^\perp c d x \\
&=& p^\perp (a + b + cd) x + x (a + b + cd) p, 
 \end{eqnarray*} 
giving the result.   

For the direct statement, we 
will use the notation of Section \ref{Section:Glimm-Halpern} without further explanation. By Lemma \ref{Chatterjee-Smith},
\[
	a \otimes 1 + 1 \otimes b + c \otimes d = \sum_{k=1}^n (x_k z_k \otimes
		y_k - x_k \otimes z_k y_k)
\]
for some $x_1, x_2, ..., x_n, y_1, y_2, ..., y_n \in \mathcal{R}$ and some $z_1, z_2, ..., z_n \in \mathcal{Z(R)}$. Fix $\omega \in \Omega$ and apply the map $\mathcal{R} \otimes \mathcal{R} \to \mathcal{R}_\omega \otimes \mathcal{R}_\omega: x \otimes y \mapsto x(\omega) \otimes y(\omega)$ to the previous equation. One obtains that
\begin{eqnarray*}
	a(\omega) \otimes 1(\omega) &+& 1(\omega) \otimes b(\omega) + 
		c(\omega) \otimes d(\omega)\\
	&=& \sum_{k=1}^n [x_k(\omega)z_k(\omega) \otimes y_k(\omega) -
		x_k(\omega) \otimes z_k(\omega)y_k(\omega)]\\
	&=& 0
\end{eqnarray*}
by Corollary \ref{Glimm_cor_2}.
If $1(\omega) \neq 0$, Lemma \ref{vector_sp_lem} tells us that either $a(\omega)$ and $c(\omega)$ are scalars, or $b(\omega)$ and $d(\omega)$ are scalars. If $1(\omega) = 0$, then $\mathcal{R}_\omega = 0$ and so $a(\omega) = b(\omega) = c(\omega) = d(\omega) = 0$. For each $x \in \mathcal{R}$, set
\[
	F_x = \{\omega \in \Omega: x(\omega) \text{ is a scalar}\}.
\]
Since $\mathcal{R}_\omega$ has a faithful irreducible representation for all $\omega \in \Omega$ (Theorem \ref{Halpern_thm_1}), $\mathcal{Z}(\mathcal{R}_\omega) = \mathbb{C}$ for all $\omega \in \Omega$. Hence
\begin{eqnarray*}
	F_x &=& \{\omega \in \Omega: x(\omega)y(\omega) = y(\omega)x(\omega)
		\text{ for all } y \in \mathcal{R}\}\\
	&=& \{\omega \in \Omega: \|x(\omega)y(\omega) - y(\omega)x(\omega)\|
		= 0 \text{ for all } y \in \mathcal{R}\}\\
	&=& \{\omega \in \Omega: \|(xy - yx)(\omega)\| = 0 \text{ for all }
		y \in \mathcal{R}\}
\end{eqnarray*}
is closed by Theorem \ref{Glimm_thm_2}. Set $F = F_a \cap F_c$. Then $G = \Omega \backslash F$ is open $\Rightarrow$ $\overline{G}$ is both open and closed. We have that $G \subset F_b \cap F_d$ $\Rightarrow$ $\overline{G} \subset F_b \cap F_d$. Also, $\Omega \backslash \overline{G} \subset \Omega \backslash G = F = F_a \cap F_c$. Let $f = 1 - \chi_{\overline{G}}$. Then $f \in C(\Omega)$. Thus (Corollary \ref{Glimm_cor_3}) there exists a $p \in \mathcal{Z(R)}$ such that $p(\omega) = f(\omega)$ for all $\omega \in \Omega$. Clearly $p$ is a projection. Since
\[
	(pa)(\omega) = p(\omega)a(\omega) = f(\omega)a(\omega) \in \mathbb{C}
\]
for all $\omega \in \Omega$, $pa \in \mathcal{Z(R)}$ (Corollary \ref{Glimm_cor_2}). Likewise, $pc, p^\perp b, p^\perp d \in \mathcal{Z(R)}$.
\end{proof}

A shorter proof of the last theorem, avoiding reduction
theory, may be given using the
Dixmier approximation theorem \cite{KR_II}.  However the 
above arguments illustrate the Glimm-Halpern techniques, which we feel will
be useful elsewhere.

\section{One-Sided $M$-Projections on $C^*$-Algebras}

In this section, we give all but one of the  promised 
%`sharper'
characterizations of complete one-sided $M$-projections on $C^*$-algebras.
The following result, which gives the 
equivalence of (i) and (iii) in Theorem 
\ref{main},  was used in \cite{BEZ}.   In the statement below,
$LM(\mathcal{A})$ is the left multiplier algebra \cite{Ped}. 

\begin{lemma}
\label{vonNeum} 
Suppose that $\mathcal{A}$ is a $C^*$-algebra. Then a linear mapping $
\varphi : \mathcal{A} \rightarrow 
\mathcal{A}$ has the form $\varphi (x)=bx$ for some $b\in LM(\mathcal{A})$
with $\left\| b\right\| \leq 1$ if and only if the column mapping
\begin{equation*}
\tau _{\varphi }^{c}:C_{2}(\mathcal{A})\rightarrow C_{2}(\mathcal{A}):
\begin{bmatrix}
x \\
y
\end{bmatrix}
\mapsto
\begin{bmatrix}
\varphi (x) \\
y
\end{bmatrix}
\end{equation*}
is contractive. Moreover if these conditions hold, then $b$ is an orthogonal 
projection in $M(\mathcal{A})$ if and only if $\varphi$ is an idempotent 
linear map. 
\end{lemma}
 
\proof    Suppose that  $\varphi (x)=bx$ for some $b\in  LM(\mathcal{A})$
with $\left\| b\right\| \leq 1$.   The fact that 
$x^*b^* bx + y^* y \leq x^* x + y^* y$, for any $x, y \in \mathcal{A}$,
yields immediately that $\tau _{\varphi }^{c}$ is contractive. 
It is fairly clear that  $b^2 = b$ if and only if
this $\varphi$ is idempotent, and since $\left\| b\right\| \leq 1$
it follows that $b^* = b \in M(\mathcal{A})$.

For the remaining implication, we suppose that
$\tau _{\varphi }^{c}$ is contractive.
We first claim that it 
suffices to prove the result in the case that $\mathcal{A}$ is a von 
Neumann algebra.  To see this, 
note that the second 
dual $(\tau _{\varphi }^{c})^{**}$ is contractive,
and may be equated with the map
$\tau _{\varphi^{**}}^{c}$ on $C_2(\mathcal{A}^{**})$.   If the 
lemma is true in the 
von Neumann algebra case, then  
there exists an element
 $b \in \mathcal{A}^{**}$ with $\|b\| \leq 1$ such that 
$\varphi^{**}(x) = b x$ for all $x \in \mathcal{A}^{**}$.
Thus $\varphi(a) = b a$ for $a \in \mathcal{A}$. 
Since $\varphi(a) \in \mathcal{A}$, we see that $b \in LM(\mathcal{A})$.
 
Henceforth assume that $\mathcal{A}$ is a von
Neumann algebra, and apply
$\tau _{\varphi }^{c}$ to the column in $C_2(\mathcal{A})$ with entries
$e$ and $1-e$, for an orthogonal
projection $e \in \mathcal{A}$.   We obtain
$\varphi (e)^{*}\varphi (e)+(1-e)\leq 1$ and thus
\begin{equation*}
(1-e)\varphi (e)^{*}\varphi (e)(1-e)=0,
\end{equation*} giving
$\varphi (e)(1-e)=0.$ But this
relation also holds for the projection $1-e,$
i.e., we have $\varphi (1-e)e=0.$ We conclude that
\begin{equation*}
\varphi (e)=\varphi (e)e=\varphi (1)e.
\end{equation*}
Since the linear span of the
projections is norm dense in $\mathcal{A},$ $\varphi
(x)=bx$ for all $x\in \mathcal{A}$,
where $b = \varphi (1)$.
\endproof

As we pointed out in \cite{BEZ}, unitary elements in 
a $C^*$-algebra may be characterized in terms
of the maps $\varphi$ on $\mathcal{A}$ such that  $\tau_\varphi$ 
is a surjective isometry.

The following result, a `sharpening' of the `order-bounded'
condition in Section 4 of \cite{Ble_Shilov}, gives the 
equivalence of (i) and (v) in Theorem \ref{main}:
 
\begin{corollary}  \label{bshc}   Suppose that $\mathcal{A}$ is a 
$C^*$-algebra. Then a linear mapping $
\varphi : \mathcal{A} \rightarrow \mathcal{A}$ 
has the form $\varphi (x)=bx$ for some $b\in LM(\mathcal{A})$
with $\left\| b\right\| \leq 1$ if and only if
for all $x \in \mathcal{A}$, there is a 2-isometric
linear embedding  $\sigma : \mathcal{A} \rightarrow B(\mathfrak{H})$ such that
$\sigma(\varphi(x))^* \sigma(\varphi(x)) \leq \sigma(x)^* \sigma(x)$.
This is also equivalent to saying that 
$\varphi(x)^* \varphi(x)  \leq x^* x$ for all $x \in \mathcal{A}$.   
\end{corollary}  

\begin{proof}   If $\varphi (x)=bx$ for some $b\in LM(\mathcal{A})$
with $\left\| b\right\| \leq 1$,
then $\varphi(x)^*\varphi(x) = x^* b^* b x \leq x^* x$ for all $x \in \mathcal{A}$.

Conversely, if
 $x \in \mathcal{A}$, and if a 2-isometric $\sigma$ exists as above,
then for $y \in \mathcal{A}$ we have
\begin{eqnarray*}
\left\| \begin{bmatrix}
                        \varphi(x)\\
                        y   
                \end{bmatrix} \right\|^2 
 &=&  \left\| \begin{bmatrix} \sigma(\varphi(x)) \\ \sigma(y)
 \end{bmatrix} \right\|^2  \\
&=&  \| \sigma(\varphi(x))^* \sigma(\varphi(x)) 
+ \sigma(y)^* \sigma(y) \| \\
 & \leq & \| \sigma(x)^* \sigma(x) + \sigma(y)^* \sigma(y) \| \\ 
&=& \left\| \begin{bmatrix} x \\  y
                \end{bmatrix} \right\|^2 .  
\end{eqnarray*}
By the previous result, $\varphi (x)=bx$ for some $b\in LM(\mathcal{A})$ with $\|b\| \leq 1$.  
\end{proof}

Theorem \ref{Smith} will give our next  
characterization, the equivalence of (i) and (ii) in Theorem 
\ref{main}. First we need some preliminary results. 

We recall that
an element $b$ of a unital 
Banach algebra $\mathcal{B}
$ is said to be $\emph{Hermitian}$ if $\|e^{itb}\| = 1$ for all $t \in \mathbb{R}$. Accordingly, a bounded linear map $T$ of a Banach space $X$ into itself is said to be Hermitian if it is a Hermitian element of the unital Banach algebra $B(X)$.
 
\begin{lemma}[Sakai-Sinclair] \label{Sakai-Sinclair}
Let $\mathcal{R}$ be a von Neumann algebra. If $T:\mathcal{R} \to \mathcal{R}$ is Hermitian, then there exist $h, k \in \mathcal{R}_{sa}$ such that $Tx = hx + xk$ for all $x \in \mathcal{R}$.
\end{lemma}
 
\begin{proof}
By \cite{Sin70}, there exists an $a \in \mathcal{R}_{sa}$ and a *-derivation $D:\mathcal{R} \to \mathcal{R}$ such that $Tx = ax + D(x)$ for all $x \in \mathcal{R}$. For the reader's convenience, we recall that a derivation is a linear map such that $D(xy) = xD(y) + D(x)y$ for all $x, y \in \mathcal{R}$. A *-derivation
satisfies the additional requirement that $D(x^*) = -D(x)^*$ for all $x \in \mathcal{R}$. By \cite{Sak66}, $D$ is inner. That is, there exists a $b \in \mathcal{R}$ such that $D(x) = bx - xb$ for all $x \in \mathcal{R}$. We claim that $b -
b^* \in \mathcal{Z(R)}$. Indeed, for all $x \in \mathcal{R}$,
\[
        bx^* - x^*b = D(x^*) = -D(x)^* = -(bx - xb)^* = b^*x^* - x^*b^*
\]
and so
\[
        (b - b^*)x^* = x^*(b - b^*).
\]
Now let $b = c + id$, with $c, d \in \mathcal{R}_{sa}$. Then it is easy to verify that $d \in \mathcal{Z(R)}$. Thus,
\[
        Tx = ax + D(x) = ax + (c + id)x - x(c + id) = (a + c)x - xc = hx + xk
\]
for all $x \in \mathcal{R}$, where $h = a + c$ and $k = -c$ are elements of $\mathcal{R}_{sa}$.
\end{proof}

\begin{lemma} \label{Hermitian}
Let $X$ be an operator space and let $P:X \to X$ be a linear idempotent. If 
$\nu_P^c:X \to C_2(X)$ is an isometry, then $P$ is Hermitian.
\end{lemma}

\begin{proof}
Let $t \in \mathbb{R}$. Then for any $x \in X$,
\begin{eqnarray*}
	\|e^{itP}x\| &=& \left\|\sum_{k=0}^\infty \frac{(itP)^k}{k!}x\right\|
		= \left\|x + \sum_{k=1}^\infty \frac{(it)^k}{k!}Px\right\|\\
	&=& \|x + (e^{it} - 1)Px\| = \|\nu_P^c(x + 
		(e^{it} - 1)Px)\|_{C_2(X)}\\
	&=& \left\|
		\begin{bmatrix}
			Px + (e^{it} - 1)Px\\
			(Id - P)x
		\end{bmatrix}
		\right\|_{C_2(X)} = \left\|
		\begin{bmatrix}
			e^{it}Px\\
			(Id - P)x
		\end{bmatrix}
		\right\|_{C_2(X)}\\
	&=& \left\|
		\begin{bmatrix}
			e^{it} & 0\\
			0 & 1
		\end{bmatrix}
		\begin{bmatrix}
			Px\\
			(Id - P)x
		\end{bmatrix}
		\right\|_{C_2(X)} = \left\|
		\begin{bmatrix}
			Px\\
			(Id - P)x
		\end{bmatrix}
		\right\|_{C_2(X)}\\
	&=& \|\nu_P^c(x)\|_{C_2(X)} = \|x\|.
\end{eqnarray*}
Hence $\|e^{itP}\| = 1$. Since the choice of $t$ was arbitrary, $P$ is Hermitian.
\end{proof}

In the language of \cite{BEZ}, the last result says that 
any `one-sided $M$-projection' (as opposed to `complete one-sided $M$-projection') is Hermitian.   An almost identical argument
shows that any `one-sided $L$-projection' is Hermitian.

\begin{theorem} \label{characterization1}
Let $\mathcal{A}$ be a $C^*$-algebra and let $P:\mathcal{A} \to \mathcal{A}$ be a linear idempotent. Then $P$ is a complete left $M$-projection if and only if $\nu_P^c:\mathcal{A} \to C_2(\mathcal{A})$ is an isometry. 
\end{theorem}

\begin{proof}
If $P$ is a complete left $M$-projection, then $\nu_P^c$ is a complete isometry. For the converse, it suffices   
(as in Lemma  
\ref{vonNeum}) to prove the 
assertion under the assumption that $\mathcal{A} = \mathcal{R}$, a von Neumann algebra. Indeed, if the assertion is true for von Neumann algebras then given a $C^*$-algebra $\mathcal{A}$ and a linear idempotent $P:\mathcal{A} \to \mathcal{A}$ such that $\nu_P^c:\mathcal{A} \to C_2(\mathcal{A})$ is an isometry, 
we have that $(\nu_P^c)^{**}:\mathcal{A}^{**} \to C_2(\mathcal{A})^{**}$ is an isometry.
This implies that
 $\nu_{P^{**}}^c:\mathcal{A}^{**} \to C_2(\mathcal{A}^{**})$ is an isometry,
so that $P^{**}$ is a complete left $M$-projection.   Thus
$P$ is a complete left $M$-projection by \cite{BEZ}, Corollary 3.5,
for example. 

By Lemma \ref{Hermitian} and Theorem \ref{Sakai-Sinclair}, there exist $h, k \in \mathcal{R}_{sa}$ such that $Px = hx + xk$ for all $x \in \mathcal{R}$. Since $P$ is idempotent,
\[	
	(h^2 - h)x + x(k^2 - k) + 2hxk = 0
\]
for all $x \in \mathcal{R}$. By Theorem \ref{Smith}, there exists a central projection $p \in \mathcal{R}$ such that $p(h^2 - h), ph, p^\perp (k^2 - k), p^\perp k \in \mathcal{Z(R)}$. Thus
\begin{eqnarray*}
	Px &=& phx + p^\perp hx + xpk + xp^\perp k\\
	&=& p^\perp (h + k)x + xp(h + k)\\
	&=& p^\perp ex + xpe,
\end{eqnarray*}
where $e = h + k \in \mathcal{R}_{sa}$. Note that $P(1) = e$, so that
\begin{eqnarray*}
	e &=& P(1) = P^2(1) = P(e) = P(h + k)\\
	&=& h^2 + 2hk + k^2 = (h + k)^2 + (hk - kh)\\
	&=& e^2 + (hk - kh).
\end{eqnarray*}
Since $e$ is self-adjoint, we also have that $e = e^2 + (kh - hk)$. Averaging, we obtain that $e = e^2$ and so $e$ is a projection. Also, 
\begin{eqnarray*}
	\|x\|^2 &=& \|\nu_P^c(x)\|_{C_2(\mathcal{R})}^2 = \left\|
		\begin{bmatrix}
			Px\\
			(Id - P)x
		\end{bmatrix}
		\right\|_{C_2(\mathcal{R})}^2\\
	&=& \|[Px]^*[Px] + [(Id - P)x]^*[(Id - P)x]\|\\
	&=& \|(p^\perp ex + xpe)^*(p^\perp ex + xpe) +
		(p^\perp e^\perp x + xpe^\perp)^*(p^\perp e^\perp x 
		+ xpe^\perp)\|\\
	&=& \|x^*exp^\perp + ex^*xep + x^*e^\perp xp^\perp + 
		e^\perp x^*xe^\perp p\|\\
	&=& \|x^*xp^\perp + (ex^*xe + e^\perp x^*xe^\perp)p\|\\
	&=& \max\{\|xp^\perp\|,\|xep\|,\|xe^\perp p\|\}^2.
\end{eqnarray*}
for all $x \in \mathcal{R}$. In particular, for all $x \in \mathcal{R}p$,
\[
	\|x\| = \max\{\|xe\|,\|xe^\perp\|\},
\]
which says that the map $\mathcal{R}p \to \mathcal{R}p:x \mapsto xe$ is a (classical) $M$-projection. Thus $ep \in \mathcal{Z}(\mathcal{R}p)$ (cf. V.4.7 in \cite{HWW93}). But then
\[
	Px = p^\perp ex + xpe = p^\perp ex + pex = ex
\]
for all $x \in \mathcal{R}$. Hence, $P$ is a complete left $M$-projection.
\end{proof}

\begin{remark}
This theorem is not valid for general operator spaces. For any Hilbert space $\mathfrak{H}$ and any nontrivial projection $P \in B(\mathfrak{H})$, $P$ (regarded as a linear map of the operator space $X = \max(\mathfrak{H})$ into itself) is an idempotent such that $\nu_P^c:X \to C_2(X)$ is an isometry. Yet $P$ is not a complete left $M$-projection (\cite{BEZ}, Proposition 6.10).
The same example shows that Lemma \ref{vonNeum} 
is not valid for general operator spaces. 
\end{remark}

The equivalence of (i) and (iv) in Theorem \ref{main} is taken up in Section 6. The equivalence of (i) and (vi) is straightforward and is left to the interested reader (hint: (i) implies (vi) is trivial and (vi) implies (v) by the `Pythagorean theorem').

\section{One-Sided $M$-Projections in von Neumann Algebra Preduals}

Using the techniques of the previous section, we can 
prove that von Neumann algebra preduals and duals have trivial complete one-sided $M$-structure.
 
\begin{theorem} \label{predual}
Let $\mathcal{R}$ be a von Neumann algebra. Then $\mathcal{R}$ has no nontrivial complete right $L$-projections. Hence, $\mathcal{R}_*$ has no nontrivial complete right $M$-ideals.
\end{theorem}

\begin{proof}
First we prove the second assertion. So assume the first assertion and suppose $J \subset \mathcal{R}_*$ is a nontrivial complete right $M$-ideal. Then $J^\perp \subset \mathcal{R}$ is a complete left $L$-summand (\cite{BEZ}, Corollary 3.6). Let $P:\mathcal{R} \to \mathcal{R}$ be the corresponding complete right $L$-projection. By assumption, $P = 0 \text{ or } Id$. Thus, $J^\perp = \{0\} \text{ or } \mathcal{R}$ $\Rightarrow$ $J = \overline{J}^{\wk} = J^{\perp\perp} \cap \mathcal{R}_* = \{0\} \text{ or } \mathcal{R}_*$, a contradiction.

Now we prove the first assertion. So suppose $P:\mathcal{R} \to \mathcal{R}$ is a complete right $L$-projection. Then $P$ is Hermitian
(by the remark after Lemma \ref{Hermitian}; or by  
Lemma \ref{Hermitian} and the duality results in \cite{BEZ}). 
Arguing exactly as in the proof of Theorem \ref{characterization1}, there exist projections $e, p \in \mathcal{R}$, with $p$ central, such that
\[
	Px = p^\perp ex + xpe
\]
for all $x \in \mathcal{R}$. Since $P$ is nontrivial, $e$ is nontrivial. Thus there exist unit vectors $\xi, \eta \in \mathfrak{H}$ (the Hilbert space on which $\mathcal{R}$ acts) such that $e\xi = \xi$ and $e\eta = 0$. Define $f \in \mathcal{R}^*$ by
\[
	f(x) = \langle x\xi,\xi \rangle + \langle x\eta,\eta \rangle
\]
for all $x \in \mathcal{R}$. Then
\[
	\|f\| \geq |f(1)| = |\langle \xi,\xi \rangle + 
		\langle \eta,\eta \rangle| = 2.
\]
On the other hand,
\begin{eqnarray*}
	|P^*(f)(x)| &=& |f(Px)| = |f(p^\perp ex) + f(xpe)|\\
	&=& |\langle p^\perp ex\xi,\xi \rangle + 
		\langle p^\perp ex\eta,\eta \rangle + \langle xpe\xi,\xi 
		\rangle + \langle xpe\eta,\eta \rangle|\\
	&=& |\langle p^\perp x\xi,\xi \rangle + \langle xp\xi,\xi \rangle| =
		|\langle x\xi,\xi \rangle| \leq \|x\|
\end{eqnarray*}
and
\begin{eqnarray*}
	|(Id - P^*)(f)(x)| &=& |f((Id - P)x)|\\
	&=& |f(p^\perp e^\perp x) + f(xpe^\perp)|\\
	&=& |\langle p^\perp e^\perp x\xi,\xi \rangle + 
		\langle p^\perp e^\perp x\eta,\eta \rangle +
		\langle xpe^\perp \xi,\xi \rangle + 
		\langle xpe^\perp\eta,\eta \rangle|\\
	&=& |\langle p^\perp x\eta,\eta \rangle + \langle xp\eta,\eta \rangle|
		= |\langle x\eta,\eta \rangle| \leq \|x\|
\end{eqnarray*}
for all $x \in \mathcal{R}$.   
Thus
 $\|P^*(f)\|, \|(Id - P^*)f\| \leq 1$, which leads to the contradiction
\begin{eqnarray*}
	2 &\leq& \|f\| = \|\nu_{P^*}^c(f)\|_{C_2(\mathcal{R}^*)} = \left\|
		\begin{bmatrix}
			P^*f\\
			(Id - P^*)f
		\end{bmatrix}
		\right\|_{C_2(\mathcal{R}^*)}\\ 
	&\leq& \sqrt{\|P^*f\|^2 + \|(Id - P^*)f\|^2} \leq \sqrt{2}.
\end{eqnarray*}
\end{proof}

We remark that an analysis of the proof shows in fact 
that there are no nontrivial `strong right $L$-projections',
in the language of \cite{BEZ}, on a
von Neumann algebra.
Indeed 
there are no nontrivial projections $P$ on a
von Neumann
algebra, such that $P^*$ is a one-sided $M$-projection. 

The following result, the classical analog of Theorem \ref{predual}, is surely well-known. However, we could not find a reference for it. While it may be proven by a slight modification to the proof of Theorem \ref{predual}, we will give an elementary argument. Note the crucial use of extreme points of classical $L$-summands, a technology which is not available yet in the one-sided theory. 

\begin{proposition}
Let $\mathcal{R}$ be a von Neumann algebra. Then $\mathcal{R}$ has 
no nontrivial classical $L$-summands.
\end{proposition}
 
\begin{proof}
Suppose that $\mathcal{R} = X \oplus_1 Y$. By Lemma 1.5 in \cite{HWW93},
\[
        \ext(\ball(\mathcal{R})) = \ext(\ball(X)) \cup \ext(\ball(Y)).
\]
Thus for every $u \in \mathcal{U(R)}$, the unitary group of $\mathcal{R}$, either $u \in X$ or $u \in Y$. Without loss of generality, we may assume $1 \in Y$. Now let $u \in \mathcal{U(R)} \cap X$. For any $\theta \in \mathbb{R}$, one has that
\[
        \|u + e^{i\theta}1\| = \|u\| + \|e^{i\theta}1\| = 2.
\]
On the other hand,
\[
        \|u + e^{i\theta}1\| = r(u + e^{i\theta}1) = \sup_{\lambda \in
                \sigma(u)} |\lambda + e^{i\theta}|.
\]
Thus, we must have that $e^{i\theta} \in \sigma(u)$. Since the choice of $\theta$ was arbitrary, $\sigma(u) = \mathbb{T}$. Now let $a \in \ball(\mathcal{R}_{sa})$. Then $a = \frac{u + u^*}{2} = \Real(u)$, where $u = a + i\sqrt{1 - a^2} \in
\mathcal{U(R)}$. Suppose $u \in X$. Then, as we just saw, $\sigma(u) =
\mathbb{T}$, so that $\sigma(a) = \sigma(\Real(u)) = \Real(\sigma(u)) = 
[-1,1]$. Likewise, since $a = \Real(u^*)$, one has that if 
$u^* \in X$ then $\sigma(a) = [-1,1]$. 
Thus if  $a$ is a projection in $\mathcal{R}$, 
it must be that $a \in Y$.  Consequently
$Y = \mathcal{R}$, and $X = \{0\}$.
\end{proof}

As a consequence of the previous two results, we have:

\begin{corollary}
Let $\mathcal{A}$ be a $C^*$-algebra. Then $\mathcal{A}$ has no nontrivial classical $L$-projections or complete right $L$-projections.
\end{corollary}

\begin{proof}
Let $P:\mathcal{A} \to \mathcal{A}$ be a complete right $L$-projection. Then $P^{**}:\mathcal{A}^{**} \to \mathcal{A}^{**}$ is also a complete right $L$-projection (\cite{BEZ}, Corollary 3.5). By Theorem \ref{predual}, $P^{**} = 0 \text{ or } Id$. Thus $P = P^{**}|_{\mathcal{A}} = 0 \text{ or } Id$. A similar proof gives the other assertion.
\end{proof}

In fact (cf. the discussion after Theorem \ref{predual}), there are no nontrivial projections $P$ on a $C^*$-algebra such that $P^*$ is a one-sided $M$-projection.

\section{One-Sided Pseudo-Orthogonality}

In this section we give our final characterization of the complete one-sided 
$M$-projections on a $C^*$-algebra (the equivalence of (i) and (iv) in
Theorem \ref{main}). This characterization is also valid for 
operator algebras having a contractive approximate identity. 
The definitions and results of this section, with the exception of the 
last two corollaries,  
are from the thesis of the last author (\cite{Zar01}).\\

Let $X$ be an operator space and $x, y \in X$. We recall
\cite{BEZ} that
 $x$ and $y$ are \emph{left orthogonal} (written $x \perp_L y$) if there exists a complete isometry $\sigma:X \to B(\mathfrak{H})$ such that $\sigma(x)^*\sigma(y) = 0$. We say that $x$ and $y$ are \emph{left pseudo-orthogonal} (and write $x \top_L y$)
 if $\left\|\begin{bmatrix} x & y \end{bmatrix}\right\|_{R_2(X)} = \max\{\|x\|,\|y\|\}$. We have that $x \perp_L y \Rightarrow x \top_L y$. Indeed, if $\sigma:X \to B(\mathfrak{H})$ is a complete isometry 
(or even a 2-isometry) such that $\sigma(x)^*\sigma(y) = 0$, then
\begin{eqnarray*}
	\left\|\begin{bmatrix} x & y \end{bmatrix}\right\|_{R_2(X)}^2 &=&
		\left\|\begin{bmatrix} \sigma(x) & \sigma(y) 
		\end{bmatrix}\right\|_{R_2(B(\mathfrak{H}))}^2\\
	&=& \left\|
		\begin{bmatrix}
			\sigma(x)^*\\
			\sigma(y)^*
		\end{bmatrix}
		\begin{bmatrix}
			\sigma(x) & \sigma(y)
		\end{bmatrix}
		\right\|_{M_2(B(\mathfrak{H}))}\\
	&=& \left\|
		\begin{bmatrix}
			\sigma(x)^*\sigma(x) & \sigma(x)^*\sigma(y)\\
			\sigma(y)^*\sigma(x) & \sigma(y)^*\sigma(y)
		\end{bmatrix}
		\right\|_{M_2(B(\mathfrak{H}))}\\
	&=& \left\|
		\begin{bmatrix}
			\sigma(x)^*\sigma(x) & 0\\
			0 & \sigma(y)^*\sigma(y)
		\end{bmatrix}
		\right\|_{M_2(B(\mathfrak{H}))}\\
	&=& \max\{\|\sigma(x)\|,\|\sigma(y)\|\}^2 = \max\{\|x\|,\|y\|\}^2.
\end{eqnarray*}
In \cite{BEZ}, Theorem 5.1, it is shown that an idempotent linear map $P:X \to X$ is a complete left $M$-projection if and only if $Px \perp_L (Id - P)y$ for all $x, y \in X$. Combining this with the previous observation, if $P$ is a complete left $M$-projection then $Px \top_L (Id - P)y$ for all $x, y \in X$. The converse is false in general (Remark \ref{counter}), but true for operator algebras having a contractive approximate identity (Theorem \ref{Zarikian}).\\

We begin with some lemmas. The first lemma concerns the concept of ``peaking''. We say that an operator $x \in B(\mathfrak{H})$ \emph{peaks} at $\xi \in \mathfrak{H}$ if $\|\xi\| = 1$ and $\|x\xi\| = \|x\|$.

\begin{lemma} \label{VZ1}
Let $x, y \in B(\mathfrak{H})$. If $x \top_L y$, $\|x\| \geq \|y\|$, and $x$ peaks at $\xi \in \mathfrak{H}$, then $x\xi \perp y\mathfrak{H}$.
\end{lemma}

\begin{proof} We may assume that $x \neq 0$, for otherwise the lemma is trivially true. Suppose that it is false that $x\xi \perp y\mathfrak{H}$. Then there exists a unit vector $\eta \in \mathfrak{H}$ such that $y\eta = \alpha x\xi + \zeta$, where $\alpha > 0$ and $\zeta \in \mathfrak{H}$ is orthogonal to $x\xi$. For all $\beta, \gamma > 0$, we have that
\begin{eqnarray*}
	\left\|
	\begin{bmatrix}
		x & y
	\end{bmatrix}
	\begin{bmatrix}
		\beta\xi\\
		\gamma\eta
	\end{bmatrix}
	\right\|^2 &=& \|\beta x\xi + \gamma y\eta\|^2\\ &=&
	\|\beta x\xi + \gamma(\alpha x\xi + \zeta)\|^2\\ &=&
	\|(\beta + \gamma\alpha)x\xi + \gamma\zeta\|^2\\ &=&
	(\beta + \gamma\alpha)^2\|x\xi\|^2 + \gamma^2\|\zeta\|^2\\ &\geq&
	(\beta + \gamma\alpha)^2\|x\|^2
\end{eqnarray*}
and
\[
	\left\|
	\begin{bmatrix}
		\beta\xi\\
		\gamma\eta
	\end{bmatrix}
	\right\|^2 = \beta^2\|\xi\|^2 + \gamma^2\|\eta\|^2 = \beta^2 + \gamma^2
\]
Thus,
\[
	\left\|
	\begin{bmatrix}
		x & y
	\end{bmatrix}
	\right\|_{R_2(B(\mathfrak{H}))}^2 \geq \frac{(\beta + 
		\gamma\alpha)^2\|x\|^2}{\beta^2 + \gamma^2}
\]
for all $\beta, \gamma > 0$. Choosing $\beta = 1/\alpha$ and $\gamma = 1$ yields
\[
	\left\|
	\begin{bmatrix}
		x & y
	\end{bmatrix}
	\right\|_{R_2(B(\mathfrak{H}))}^2 \geq (1 + \alpha^2)\|x\|^2 > \|x\|^2 = 	\max\{\|x\|,\|y\|\}^2,
\]
a contradiction. Thus $x\xi \perp y\mathfrak{H}$. \end{proof}

\begin{lemma} \label{VZ2}
Let $X \subset B(\mathfrak{H})$ be an operator space and let $P:X \to X$ be an idempotent linear map such that $Px \top_L (Id - P)y$ for all $x, y \in X$. Set $Y = PX$ and $Z = (Id - P)X$. If $y \in Y$ peaks at $\xi \in \mathfrak{H}$, then $y\xi \perp Z\mathfrak{H}$.
\end{lemma}

\begin{proof}
We may assume that $y \neq 0$, for otherwise the lemma is trivially true. Let $z \in Z$. There exists a $\kappa > 0$ such that $\|\kappa y\| \geq \|z\|$. Since $\kappa y \in Y$, $\kappa y \top_L z$. Since $y$ peaks at $\xi$, the same is true for $\kappa y$. By Lemma \ref{VZ1}, $\kappa y\xi \perp z\mathfrak{H}$. Thus, $y\xi \perp z\mathfrak{H}$. Since the choice of $z$ was arbitrary, $y\xi \perp Z \mathfrak{H}$. \end{proof}

\begin{lemma} \label{VZ3}
Let $X$ be an operator space and let $P:X \to X$ be an idempotent linear map such that $Px \top_L (Id - P)y$ for all $x, y \in X$. Then $P$ is contractive.
\end{lemma}

\begin{proof}
By Ruan's theorem (\cite{ER_book}, Theorem 2.3.5), we may assume that $X \subset B(\mathfrak{H})$. By passing to the universal representation of $B(\mathfrak{H})$, if necessary, we may assume that every $x \in X$ peaks at some $\xi \in \mathfrak{H}$. Now let $x \in X$. Then $x = y + z$, where $y = Px$ and $z = (Id - P)x$. Let $\xi \in \mathfrak{H}$ be such that $y$ peaks at $\xi$. Then $y\xi \perp z\mathfrak{H}$ by Lemma \ref{VZ2}. Hence,
\[
	\|Px\| = \|y\| = \|y\xi\| \leq \sqrt{\|y\xi\|^2 + \|z\xi\|^2} =
		\|y\xi + z\xi\| = \|x\xi\| \leq \|x\|.
\]
Since the choice of $x$ was arbitrary, $P$ is contractive. \end{proof}

\begin{theorem} \label{Zarikian}
Let $\mathcal{B}$ be a unital operator algebra and $P:\mathcal{B} \to \mathcal{B}$ be an idempotent linear map. Then $P$ is a complete left $M$-projection 
%ADDED 
(that is, $Px = e x$ for a projection $e$ of norm $1$ in
$A$) if and only if $Px \top_L (Id - P)y$ for all $x, y \in \mathcal{B}$.
\end{theorem}

\begin{proof}
We have already indicated the proof of the forward implication. For the reverse implication, represent $\mathcal{B}$ faithfully on a Hilbert space $\mathfrak{H}$ in such a way that every $x \in \mathcal{B}$ peaks at some $\xi \in \mathfrak{H}$ (e.g. use the universal representation of the C*-algebra generated by $\mathcal{B}$). Let $X = P\mathcal{B}$ and $Y = (Id - P)\mathcal{B}$, so that $x \top_L y$ for all $x \in X$ and all $y \in Y$. Then there exists a unique decomposition $I = x_0 + y_0$ with $x_0 \in X$ and $y_0 \in Y$. It suffices to prove that $x_0$ is a projection in $\mathcal{B}$, because if so, setting $e = x_0$ and invoking Lemma \ref{VZ2}, we have that $e\mathfrak{H} \perp Y\mathfrak{H}$ and $(I-e)\mathfrak{H} \perp X\mathfrak{H}$. But then $ey = 0$  for all $y \in Y$ and $(I - e)x = 0$ for all $x \in X$, and so $P(x + y) = x = e(x + y)$ for all $x \in X$ and all $y \in Y$.

By Lemma \ref{VZ3}, we know that $\|x_0\| \leq 1$ and $\|y_0\| \leq 1$. Thus, the closed linear subspaces $\mathfrak{L} = \{\xi \in \mathfrak{H}: x_0 \xi = \xi\}$ and $\mathfrak{M} = \{\eta \in \mathfrak{H}: y_0 \eta = \eta\}$ of $\mathfrak{H}$ are orthogonal, since $\mathfrak{L} \perp Y\mathfrak{H} \supset \mathfrak{M}$ by Lemma \ref{VZ2}. To complete the proof it suffices to show that $\mathfrak{N} = \mathfrak{L}^{\perp} \cap \mathfrak{M}^{\perp} = \{0\}$, because if that is the case, one quickly deduces (using the fact that $x_0 + y_0 = I$) that
\[
	x_0 =
	\begin{bmatrix}
		I & 0\\
		0 & 0
	\end{bmatrix}
\]
with respect to the decomposition $\mathfrak{H} = \mathfrak{L} \oplus \mathfrak{M}$.

We claim that if $0 \neq x \in X$ peaks at $\xi \in \mathfrak{H}$, then $x\xi \in \mathfrak{L}$. By Lemma \ref{VZ2}, $x\xi \perp Y\mathfrak{H}$. Thus with respect to the decomposition $\mathfrak{H} = \mathbb{C} x\xi \oplus (\mathbb{C} x\xi)^{\perp}$,
\[
	x_0 =
	\begin{bmatrix}
		x_{11} & x_{12}\\
		x_{21} & x_{22}
	\end{bmatrix}
	\text{ and }
	y_0 =
	\begin{bmatrix}
		0 & 0\\
		y_{21} & y_{22}
	\end{bmatrix}.
\]
Combining this with the fact that $x_0 + y_0 = I$ implies that 
\[
	x_0 =
	\begin{bmatrix} 
		I & 0\\
		x_{21} & x_{22}
	\end{bmatrix}.
\]
Thus 
\[
	x_0x\xi =
	\begin{bmatrix}
		I & 0\\
		x_{21} & x_{22}
	\end{bmatrix}
	\begin{bmatrix}
		x\xi\\
		0
	\end{bmatrix} =
	\begin{bmatrix}
		x\xi\\
		x_{21}x\xi
	\end{bmatrix}.
\]
If $x_{21}x\xi \neq 0$, then 
\[
	\|x_0x\xi\|^2 = \|x\xi\|^2 + \|x_{21}x\xi\|^2 > \|x\xi\|^2,
\]
a contradiction to the fact that $\|x_0\| \leq 1$.

Now suppose that $\mathfrak{N} \neq \{0\}$. With respect to the decomposition $\mathfrak{H} = \mathfrak{L} \oplus \mathfrak{M} \oplus \mathfrak{N}$, 
\[
	x_0 =
	\begin{bmatrix}
		I & * & *\\
		0 & 0 & 0\\
		0 & * & *
	\end{bmatrix}
	\text{ and }
	y_0 =
	\begin{bmatrix}
		0 & 0 & 0\\
		* & I & *\\
		* & 0 & *
	\end{bmatrix},
\]
where the *'s represent (possibly) nonzero bounded linear operators. Because $x_0 + y_0 = I$,
\[
	x_0 =
	\begin{bmatrix}
		I & 0 & 0\\
		0 & 0 & 0\\
		0 & 0 & z
	\end{bmatrix}
	\text{ and }
	y_0 =
	\begin{bmatrix}
		0 & 0 & 0\\
		0 & I & 0\\
		0 & 0 & I - z
	\end{bmatrix},
\]
where $z \in B(\mathfrak{N})$. It must be that
\[
	x_0 y_0 =
	\begin{bmatrix}
		0 & 0 & 0\\
		0 & 0 & 0\\
		0 & 0 & z - z^2
	\end{bmatrix} 
	\neq 0.
\]
To see this, we first note that if $z = 0$, then for any $0 \neq \xi \in \mathfrak{N}$, 
\[
	y_0\xi =
	\begin{bmatrix}
		0 & 0 & 0\\
		0 & I & 0\\
		0 & 0 & I
	\end{bmatrix}
	\begin{bmatrix}
		0\\
		0\\
		\xi
	\end{bmatrix} 
	=
	\begin{bmatrix}
		0\\
		0\\
		\xi
	\end{bmatrix} 
	= \xi,
\]
which means that $\xi \in \mathfrak{M}$, a contradiction. If $z - z^2 = 0$, then for any $\xi \in \mathfrak{N}$ with $z\xi \neq 0$,
\[
	x_0z\xi =
	\begin{bmatrix}
		I & 0 & 0\\
		0 & 0 & 0\\
		0 & 0 & z
	\end{bmatrix}
	\begin{bmatrix}
		0\\
		0\\
		z\xi
	\end{bmatrix} 
	=
	\begin{bmatrix}
		0\\
		0\\
		z^2\xi
	\end{bmatrix} 
	= z^2\xi = z\xi,
\]
which means that $z\xi \in \mathfrak{L}$, a contradiction.

Now there exists a unique decomposition $x_0y_0 = x_1 + y_1$ with $x_1 \in X$ and $y_1 \in Y$. Without loss of generality, we may assume that $x_1 \neq 0$. With respect to the decomposition $\mathfrak{H} = \mathfrak{L} \oplus \mathfrak{M} \oplus \mathfrak{N}$,
\[
	x_1 =
	\begin{bmatrix}
		* & * & *\\
		0 & 0 & 0\\
		* & * & *
	\end{bmatrix}
	\text{ and }
	y_1 =
	\begin{bmatrix}
		0 & 0 & 0\\
		* & * & *\\
		* & * & *
	\end{bmatrix}.
\]
Because $x_1 + y_1 = x_0y_0$,
\[
	x_1 =
	\begin{bmatrix}
		0 & 0 & 0\\
		0 & 0 & 0\\
		* & * & *
	\end{bmatrix}.
\]
If $x_1$ peaks at $\xi_1 \in \mathfrak{H}$, then by a previous claim we have that $x_1\xi_1 \in \mathfrak{L}$. On the other hand, we plainly have that $x_1\xi_1 \in \mathfrak{N}$, a contradiction.
\end{proof}

\begin{corollary} \label{anything}
Let $\mathcal{B}$ be an operator algebra having a contractive approximate 
identity
 and $P:\mathcal{B} \to \mathcal{B}$ be an idempotent linear map. Then $P$ is a complete left $M$-projection 
%ADDED
(that is, $Px = e x$ for a projection $e$ of norm $1$ in
$LM(A)$) if and only if $Px \top_L (Id - P)y$ for all $x, y \in \mathcal{B}$.
\end{corollary}

\begin{proof}
Again, we need only prove the reverse implication.
Recall that
$\mathcal{B}^{**}$ is a unital operator algebra with respect to the Arens product (cf. \cite{ER90}, Theorem 2.1). Since $P$ is a linear idempotent such that $Px \top_L (Id - P)y$ for all $x, y \in \mathcal{B}$, the map
\[
	\Phi:P\mathcal{B} \oplus_\infty (Id - P)\mathcal{B} \to 
		R_2(\mathcal{B}):(Px,(Id - P)y) \mapsto 
		\begin{bmatrix} Px & (Id - P)y \end{bmatrix}
\]
is an isometry. By usual duality arguments, $\Phi^{**}:(P\mathcal{B} \oplus_\infty (Id - P)\mathcal{B})^{**} \to R_2(\mathcal{B})^{**}$ is also an isometry. Now under the isometric identifications
\[
	(P\mathcal{B} \oplus_\infty (Id - P)\mathcal{B})^{**} \cong 
		(P\mathcal{B})^{**} \oplus_\infty ((Id - P)\mathcal{B})^{**},
\]
\[
	(P\mathcal{B})^{**} \cong (P\mathcal{B})^{\perp\perp} = 
		P^{**}\mathcal{B}^{**},
\]
\[
	((Id - P)\mathcal{B})^{**} \cong ((Id - P)\mathcal{B})^{\perp\perp} =
		(Id - P^{**})\mathcal{B}^{**},
\]
and
\[
	R_2(\mathcal{B})^{**} \cong R_2(\mathcal{B}^{**}),
\]
$\Phi^{**}$ corresponds to the map
\begin{eqnarray*}
	\Psi&:&P^{**}\mathcal{B}^{**} \oplus_\infty (Id - P^{**})
		\mathcal{B}^{**} \to R_2(\mathcal{B}^{**}):\\
	&& (P^{**}\overline{x}, (Id - P^{**})\overline{y}) \mapsto 
		\begin{bmatrix} P^{**}\overline{x} & (Id - P^{**})\overline{y}
		\end{bmatrix}.
\end{eqnarray*}
Thus $P^{**}\overline{x} \top_L (Id - P^{**})\overline{y}$ for all $\overline{x}, \overline{y} \in \mathcal{B}^{**}$. By Theorem \ref{Zarikian}, $P^{**}$ is a complete left $M$-projection. By \cite{BEZ}, Corollary 3.5, $P$ is a complete left $M$-projection.
\end{proof}

\begin{remark} \label{counter}
Theorem \ref{Zarikian} is not true for general operator spaces. Let $X \subset M_3$ be the linear subspace spanned by
\[
	s = 
	\begin{bmatrix}
		\sqrt{2} & 0 & 0\\
		0 & 1 & 0\\
		0 & 0 & 0
	\end{bmatrix}
	\text{ and }
	t =
	\begin{bmatrix}
		0 & 0 & 0\\
		0 & 1 & 0\\
		0 & 0 & \sqrt{2}
	\end{bmatrix}.
\]
Then $P:X \to X:\alpha s + \beta t \mapsto \alpha s$ is an idempotent linear map such that $Px \top_L (Id - P)y$ for all $x, y \in X$. However $P$ is not a complete left $M$-projection (cf. \cite{Zar01}, \S 2.9).
\end{remark}

Indeed, with a little more work,
the last example reveals that the `matricial version' of
one-sided pseudo-orthogonality does not characterize one-sided
complete $M$-projections on general operator spaces.   However 
we do have the following result:

\begin{corollary}
Let $X$ be a TRO (or, equivalently, a Hilbert $C^*$-module) and let 
$P:X \to X$ be an idempotent linear map. Then $P$ is a
complete left $M$-projection 
if and only if $P_nx \top_L (Id - P_n)y$ for all
$x, y \in M_n(X)$ and all $n \in \mathbb{N}$.   If $X$ is a
right Hilbert $C^*$-module say, then such $P$ are exactly
the adjointable projections on $X$.
\end{corollary}
 
\begin{proof}
The idea is similar to that of the proof of 
Theorem 4.6 in \cite{BEZ}, so we will simply sketch the 
argument.
Suppose that
$P_nx \top_L (Id - P_n)y$ for all $x, y \in M_n(X)$ and all $n \in
\mathbb{N}$. Then it is easy to see that $P^{**}\overline{x} \top_L (Id -
P^{**})\overline{y}$ for all $\overline{x}, \overline{y} \in M_n(X^{**})$ and
all $n \in \mathbb{N}$. Thus, by Corollary 3.5 of \cite{BEZ}, we may assume
$X$ is self-dual.  
By Lemma 4.4 of \cite{BEZ}, there exists a cardinal $J$
such that $M_J(X)$ is completely isometrically isomorphic to a von Neumann
algebra. Because $P$ is completely contractive (Lemma \ref{VZ3}), the
amplification $P_J:M_J(X) \to M_J(X)$ is well-defined. Since $P_nx \top_L (Id
- P_n)y$ for all $x, y \in M_n(X)$ and all $n \in \mathbb{N}$, one has that
$P_J\left(\begin{bmatrix} x_{i,j} \end{bmatrix}\right) \top_L (Id -
P_J)\left(\begin{bmatrix} y_{i,j} \end{bmatrix}\right)$ for all
$\begin{bmatrix} x_{i,j} \end{bmatrix}, \begin{bmatrix} y_{i,j} \end{bmatrix}
\in M_J(X)$. 
By Theorem \ref{Zarikian}, $P_J$ is a complete left $M$-projection.
Fixing $j_0 \in J$ and 
restricting $P_J$ to those elements of $M_J(X)$ with zeros everywhere except
possibly the $(j_0,j_0)$ position, we conclude that $P$ is a complete left
M-projection. The converse is trivial.
\end{proof}

We saw in Section 4 that 
the condition that $\nu_P^c$ is isometric is sufficient to 
characterize complete left $M$-projections on
$C^*$-algebras.
It would be interesting to know if this was also true for
general operator algebras.
As another corollary to the methods of this section, we do at least have the following:

\begin{corollary}   
Let $\mathcal{B}$ be an operator algebra with contractive approximate identity,
and suppose that $P$ is an idempotent linear map on $\mathcal{B}$.
Then $P$
is a complete left $M$-projection if and only if $\nu_P^c:\mathcal{B} \to
C_2(\mathcal{B})$ is 2-isometric.
\end{corollary}
 
\begin{proof}
Indeed, in the case that $\nu_P^c$ is 2-isometric, one has 
for $x, y \in \mathcal{B}$ that
\begin{eqnarray*}
        \left\|
        \begin{bmatrix}
                Px & (I   d - P)y
        \end{bmatrix}
        \right\|_{R_2(\mathcal{B})} &=&
        \left\|
        \begin{bmatrix}
                \nu_P^c(Px) & \nu_P^c((Id - P)y)
        \end{bmatrix}
        \right\|_{M_2(\mathcal{B})}\\
        &=& \left\|
        \begin{bmatrix}
                P(Px) & P((Id - P)y)\\
                (Id - P)(Px) & (Id - P)((Id - P)y)\\
        \end{bmatrix}
        \right\|_{M_2(\mathcal{B})}\\
        &=&
        \left\|
        \begin{bmatrix}
                Px & 0\\
                0 & (Id - P)y
        \end{bmatrix}
        \right\|_{M_2(\mathcal{B})}\\
        &=& \max\{\|Px\|,\|(Id - P)y\|\},
\end{eqnarray*}
which says that $Px \top_L (Id - P)y$.   The result then follows
from  Corollary \ref{anything}.    
\end{proof} 

\section{Intersections of One-Sided $M$-Summands}

In the classical theory, there is 
a well known `calculus' of $M$-summands, $L$-summands, and $M$-ideals
(see section I.1 in \cite{HWW93}).  Many of these results 
go through in the quantized, one-sided case considered in 
\cite{BEZ}.   A few do fail without extra hypotheses.  
The state of this `calculus' of complete one-sided $M$-summands
and ideals will be described in a forthcoming paper. We display
next one classical result,
namely that the intersection of two $M$-summands is again an
$M$-summand, which fails for 
general complete one-sided $M$-summands.   (We remark however that if the
complete one-sided $M$-projections
corresponding to these summands commute, then this result is
valid).   Indeed 
for a unital $C^*$-algebra $\mathcal{A}$, the 
(complete) right M-summands are exactly the principal 
right ideals $p \mathcal{A}$, for a projection $p \in \mathcal{A}$.
However for projections 
$p, q \in \mathcal{A}$, it need not be the case (unless
$\mathcal{A}$ is a von Neumann algebra) that 
$(p \mathcal{A}) \cap (q \mathcal{A}) = r \mathcal{A}$
for any projection $r \in \mathcal{A}$.  
Although this may be known, we could not find this fact in the 
literature.

The following is a modification of an example shown to us by Stephen Dilworth.
 
\begin{example} \label{A}
There exists a unital $C^*$-algebra $\mathcal{A} \subset B(\mathfrak{H})$ and projections $P, Q \in \mathcal{A}$ such that
\begin{enumerate}
\item $P \wedge Q \notin \mathcal{A}$;
\item $\mathfrak{H}_1 = \ran(P \wedge Q)$ is separable
with orthonormal basis $\{\xi_n\}$;
\item The one-dimensional subprojection of $P \wedge Q$ with 
range $\mathbb{C}\xi_n$ 
lies in $\mathcal{A}$ for all $n \in \mathbb{N}$.\end{enumerate}
\end{example}
 
\begin{proof}
Let $\mathfrak{H}_1$ be a separable Hilbert space with orthonormal basis $\{\xi_n\}$, $\mathfrak{H}_2$ be a separable Hilbert space with orthonormal basis $\{\eta_n\}$, and $\mathfrak{H} = \mathfrak{H}_1 \oplus \mathfrak{H}_2$. Define
\[
        \mathfrak{H}_3 = \overline{\spn(\{\eta_{2n-1}: n \in \mathbb{N}\})}
                \subset \mathfrak{H}_2
        \text{ and }
        \mathfrak{H}_4 = \overline{\spn(\{\eta_{2n-1} + \frac{1}{n}\eta_{2n}:
                n \in \mathbb{N}\})} \subset \mathfrak{H}_2
\]
and set
\[
        \mathfrak{H}_5 = \mathfrak{H}_1 \oplus \mathfrak{H}_3 \subset
                \mathfrak{H}
        \text{ and }
        \mathfrak{H}_6 = \mathfrak{H}_1 \oplus \mathfrak{H}_4 \subset
                \mathfrak{H}.
\]
Let $P, Q \in B(\mathfrak{H})$ be the projections onto $\mathfrak{H}_5$ and $\mathfrak{H}_6$, respectively, and
let $E_n \in B(\mathfrak{H})$ be the projection onto $\mathbb{C}\xi_n$ for all $n \in \mathbb{N}$. Denote by $\mathcal{A}$ the unital $C^*$-subalgebra of $B(\mathfrak{H})$ generated by $P$, $Q$, and the family $\{E_n: n \in \mathbb{N}\}$. We claim that $P \wedge Q$, the projection onto $\mathfrak{H}_1 = \mathfrak{H}_5 \cap \mathfrak{H}_6$, is not an element of $\mathcal{A}$. Suppose, to the contrary, that $P \wedge Q \in \mathcal{A}$. Pick $\epsilon > 0$ arbitrarily. Then there exists an operator
\[
        T = \sum_{k=1}^M \alpha_k Q^{r_k}(PQ)^{s_k}P^{t_k} +
                \sum_{l=1}^N \beta_l E_l,
\]
where $r_k, t_k \in \{0, 1\}$ and $s_k \in \mathbb{N}_0$ for all $k = 1, 2, ..., M$, such that $\|P \wedge Q - T\| < \epsilon$. It follows that for any $x \in \mathfrak{H}_1 \cap \{\xi_1, \xi_2, ..., \xi_N\}^\perp$,
\[
        \left|1 - \sum_{k=1}^M \alpha_k\right|\|x\| = \left\|x - \sum_{k=1}^M
                \alpha_k x\right\| = \|(P \wedge Q - T)x\| \leq \epsilon\|x\|.
 \]
Thus
\[
        \left|1 - \sum_{k=1}^M \alpha_k\right| \leq \epsilon.
\]
On the other hand, it is easy to check that
\[
        Q^{r_k}(PQ)^{s_k}P^{t_k} \eta_{2n-1} = \left(\frac{n^2}{n^2 + 1}\right)^
                {s_k + r_k}\left[
\eta_{2n-1} + \frac{r_k}{n}\eta_{2n}\right]
\]
$\Rightarrow$
\[
        T\eta_{2n-1} = \sum_{k=1}^M \alpha_k \left(\frac{n^2}{n^2 + 1}\right)^
                {s_k + r_k}\left[\eta_{2n-1} + \frac{r_k}{n}
\eta_{2n}\right].
\]
But then
\begin{eqnarray*}
        \epsilon &\geq& \|(P \wedge Q - T)\eta_{2n-1}\| = \|T
\eta_{2n-1}\|\\
        &=& \sqrt{\left|\sum_{k=1}^M \alpha_k\left(\frac{n^2}{n^2 + 1}\right)^
                {s_k+r_k}\right|^2 + \left|\sum_{k=1}^M \alpha_k\left(\frac
                {n^2}{n^2 + 1}\right)^{n_k+r_k}\frac{r_k}{n}\right|^2}\\
        &\to& \left|\sum_{k=1}^M \alpha_k \right| \geq 1 - \epsilon
\end{eqnarray*}
as $n \to \infty$, a contradiction.
\end{proof}

\begin{example} \label{E}
There exists a
unital $C^*$-algebra $\mathcal{A}$, 
and projections
$p, q \in \mathcal{A}$, such that
$(p \mathcal{A}) \cap (q \mathcal{A}) \neq r \mathcal{A}$
for any projection $r \in \mathcal{A}$.
\end{example}
 
\begin{proof}  
Let $\mathcal{A}$ be the $C^*$-algebra from Example \ref{A}, $J_1 = P\mathcal{A}$, and $J_2 = Q\mathcal{A}$. If $J_1 \cap J_2 = R\mathcal{A}$ for some projection $R \in B(\mathfrak{H}$), 
then necessarily $R \in \mathcal{A}$ (since $\mathcal{A}$ is unital).
Since $P R = Q R = R$, we have $R \leq P \wedge Q$.
By construction, $E_n \in \mathcal{A}$. Since $E_n 
\leq P \wedge Q$, we have that
 $PE_n = E_n$ and $QE_n = E_n$, so that
 $E_n \in J_1 \cap J_2 = R\mathcal{A}$. It follows that 
$E_n \leq R$. Since $n$ was arbitrary, we 
obtain the contradiction that $R = P \wedge Q$.
\end{proof}

\end{document}